\documentclass[11pt]{amsart}
\usepackage{amsmath,amssymb,amsfonts,amsthm,verbatim}

\DeclareMathOperator{\lcm}{lcm}
\DeclareMathOperator{\Gal}{Gal}
\DeclareMathOperator{\charp}{char}

\DeclareMathOperator{\Aut}{Aut}

\newcommand{\F}{\mathbb{F}}
\newcommand{\Z}{\mathbb{Z}}
\newcommand{\Q}{\mathbb{Q}}
\newcommand{\bC}{\mathbb{C}}

\newcommand{\col}{\,{:}\,}
\renewcommand{\hat}{\widehat}
\renewcommand{\tilde}{\widetilde}
\newcommand{\mybar}[1]{#1\llap{$\overline{\phantom{\rm#1}}$}}
\renewcommand{\bar}{\mybar}

\newcommand{\Kbar}{\bar K}

\newcommand{\Qbar}{\bar\Q}
\newcommand{\Lbar}{\bar L}
\newcommand{\Mbar}{\bar M}

\theoremstyle{plain}
\newtheorem{theorem} {Theorem}  [section]
\newtheorem{proposition} [theorem] {Proposition}
\newtheorem{corollary} [theorem] {Corollary}

\newtheorem{lemma}[theorem]{Lemma}
\newtheorem{exampleme}[theorem]{Example}

\theoremstyle{definition}
\newtheorem{definitionme}[theorem]{Definition}
\newtheorem*{definition}{Definition}

\theoremstyle{remark}
\newtheorem*{remark}{Remark}

\begin{document} 
\title{Polynomials with a common composite}

\author{Robert M.~Beals}
\address{Center for Communications Research,
         805 Bunn Drive,
         Princeton, NJ 08540-1966, USA.}
\email{beals@idaccr.org}

\author{Joseph L.~Wetherell}
\address{Center for Communications Research,
         4320 Westerra Court,
         San Diego, CA 92121-1967, USA.}
\email{jlwether@alum.mit.edu}

\author{\hbox{Michael E. Zieve}}
\address{Center for Communications Research,
         805 Bunn Drive,
         Princeton, NJ 08540-1966, USA.}
\email{zieve@math.rutgers.edu}
\urladdr{http://www.math.rutgers.edu/$\sim$zieve/}

\thanks{We thank the anonymous referee for several valuable suggestions.}

\date{\today}

\begin{abstract}
Let $f$ and $g$ be nonconstant polynomials over a field $K$.  In this paper
we study the pairs $(f,g)$ for which the intersection
$K[f]\cap K[g]$ is larger than $K$.  We describe all such pairs in case
$K$ has characteristic zero, as a consequence of classical results due to Ritt.
For fields $K$ of positive characteristic we present various results, examples,
and algorithms.
\end{abstract}

\maketitle


\section{Introduction}

Let $f_1$ and $f_2$ be nonconstant polynomials over a field $K$ of characteristic~$p\ge 0$.
In this paper we examine whether $f_1$ and $f_2$ have a \emph{common composite}, i.e.,
whether there are nonconstant $u,v\in K[x]$ such that $u(f_1(x))=v(f_2(x))$.
Any such polynomial $u(f_1(x))$ is a common composite.

It turns out that there are very precise results about common
composites whose degree is not divisible by $p$.  Namely, if $f_1$ and
$f_2$ have such a common composite, then they have a common composite
of degree $\lcm(\deg(f_1),\deg(f_2))$.  Also, under the same
hypotheses, there are $g_1,g_2,r\in K[x]$ with
$\deg(r)=\gcd(\deg(f_1),\deg(f_2))$ such that $f_1=g_1\circ r$ and
$f_2=g_2\circ r$.  Further, in
Theorem~\ref{char0} we describe all possibilities for $g_1$ and $g_2$.

For common composites of degree divisible by $p$, the situation is much more complicated.
In particular, we present counterexamples to the gcd and lcm results, as well as a
sequence of examples of pairs of bounded-degree polynomials whose least-degree common
composites have degrees growing without bound.
As a substitute, we give an algorithm which quickly
determines whether there is a common composite of degree less than any fixed bound.
Further, we prove the following result describing necessary and sufficient
criteria for two polynomials to have a common composite.
In this result, $m_i(a)$ denotes the multiplicity of $x=a$
as a root of $f_i(x)-f_i(a)$, and $\Kbar$ denotes an algebraic closure of $K$.

\begin{theorem}
\label{thmintro}
Polynomials $f_1,f_2\in K[x]\setminus K$ have a common composite if and only if
there is a nonempty finite subset $A$ of $\Kbar$ which admits
a function $\ell:A\to\Z$ such that
\begin{itemize}
\item for $a\in A$ and $i\in\{1,2\}$, \,$\ell(a)/m_i(a)$ is a positive integer; and
\item for $i\in\{1,2\}$, $a\in A$, and $b\in\Kbar$, if $f_i(a)=f_i(b)$ then $b\in A$ and
$\ell(a)/m_i(a)=\ell(b)/m_i(b)$.
\end{itemize}
\end{theorem}

We suspect that `most' pairs of polynomials $(f_1,f_2)$ have no common composite.
In characteristic zero, this follows from Ritt's 1922 results \cite{Ritt}; see also
Theorem~\ref{char0}.  However, in
positive characteristic it is difficult to produce a pair of polynomials
which one can prove do not have a common composite; in fact, it has even been
conjectured that no such polynomials exist over finite fields \cite{McConnell}.
This conjecture was disproved in the 1970's via clever examples in \cite{Astrid!},
\cite{BremnerMorton} and \cite{AlexandruPopescu}.  However, the arguments in those
papers seem to apply only to very carefully chosen polynomials.  We give some
general methods for proving two polynomials have no common composite.
A special case of our results is as follows:

\begin{theorem}
\label{thmintro2}
Suppose $f_1,f_2\in K[x]\setminus K[x^p]$ and $\alpha,\beta\in\Kbar$ satisfy
$f_i(\alpha)=f_i(\beta)$ for both $i=1$ and $i=2$.  Then $f_1$ and $f_2$ have
no common composite if either of the following hold:
\begin{itemize}
\item $m_1(\alpha)m_2(\beta)\ne m_1(\beta)m_2(\alpha)$; or
\item $f_1'(\alpha)f_2'(\beta)\ne f_1'(\beta)f_2'(\alpha)$,  and
$[K(\alpha)\col K]$ is divisible by a prime greater than \,$\max(\deg(f_1),\deg(f_2))$.
\end{itemize}
\end{theorem}

For instance, one can check that the first condition implies that $x^2+x$ and $x^3+x^2$
have no common composite over $\F_2$, and the second condition implies that
$x^4+x^3$ and $x^6+x^2+x$ have no common composite over $\F_2$.  In fact,
we expect that our most general version of the second condition will apply to `most'
pairs of polynomials over a finite field.

The existence of a common composite can be reformulated in several different ways.
It is clearly equivalent to saying the intersection of the polynomial rings
$K[f_1]$ and $K[f_2]$ is strictly bigger than $K$.
We show moreover that it is also equivalent to saying the
intersection $K(f_1) \cap K(f_2)$ is strictly bigger than $K$,
i.e., it is equivalent to $f_1$ and $f_2$ having a common rational function composite.

Finally, we mention an application of the results in this paper.
Suppose $f_1$ and $f_2$ have no common composite, and assume further that $f_1$ and $f_2$
are not both functions of any polynomial of degree more than 1.
Then the $x$-resultant of $f_1(x)-u$ and $f_2(x)-v$ is an irreducible polynomial in $K[u,v]$
which is not a factor of any nonzero `variables separated' polynomial $r(u)-s(v)$
with $r,s\in K[x]$.  We know no other way to produce such irreducibles.

Various authors have considered common composites from different perspectives, using methods
involving Riemann surfaces, power series, curves and differentials, and group theory, among others.
Of special importance is Schinzel's book~\cite{Schinzel}, which contains beautiful proofs
using (in most cases) only basic properties of polynomials.  In the first few sections of
this paper, we include new proofs of some known results. 
Also, we include multiple proofs of some results, and numerous examples illustrating the
different types of phenomena that can occur.  We hope that this will lead to future work
providing more insight into the mysteries surrounding common composites of polynomials.

We now describe the organization of this paper.
In the next section we show that if two polynomials have a common composite over an
extension of $K$, then they have a common composite over $K$.  Then in Section~\ref{sec rat}
we show that existence of a common rational function composite implies existence of
a common polynomial composite.  In Section~\ref{sec deg} we give results and examples
addressing the degrees of common composites.  In the next section we use these results to
describe all polynomials which have a common composite of degree not divisible by $\charp(K)$.
In Section~\ref{sec fib} we give an algorithm which quickly determines whether there is
a common composite of degree less than some bound.  In the final three sections we give
criteria for existence or nonexistence of common composites, and in particular
we prove generalizations of Theorems~\ref{thmintro} and~\ref{thmintro2}.


\section{Reduction to the Case of Algebraically Closed Fields}
\label{sec ac}

\begin{theorem}
\label{algclosed}
If $f_1,f_2\in K[x]$ have a common composite over the
algebraic closure $\Kbar$ of $K$, then they have a common composite over $K$.
Moreover, the minimal degree of any common composite over $K$ equals the minimal degree
of any common composite over $\Kbar$.
\end{theorem}

\begin{proof}
Let $n$ be the minimal degree of any common composite over $\Kbar$.
Then there are polynomials $g_1,g_2,h$ in $\Kbar[x]$ with
\begin{equation}\label{dependence}
h=g_1\circ f_1=g_2\circ f_2
\end{equation}
such that $h$ has degree $n$.  Let $d_i$ be the degree of $g_i$. 
Equation~(\ref{dependence}) expresses a $\Kbar$-linear
dependence of the polynomials 
\[
1, f_1, f_1^2,\ldots,f_1^{d_1},f_2,f_2^2,\ldots,f_2^{d_2}.
\]
Letting $V$ be the $K$-vector space spanned by these polynomials,
we see that the $\Kbar$-vector space $\Kbar\otimes_K V$ has
the same dimension as $V$.  Thus the polynomials are linearly independent
over $\Kbar$ if and only if they are linearly dependent over $K$.
\end{proof}

\begin{corollary}
\label{sep}
If $f_1,f_2\in K[x]\setminus K[x^p]$ have a common composite, then they have a common
composite which is not in $K[x^p]$.
\end{corollary}

\begin{proof}
Let $h\in K[x]$ be a minimal degree common composite of $f_1$ and $f_2$, and assume $h\in K[x^p]$.
Write $h=g_1\circ f_1=g_2\circ f_2$ with $g_1,g_2\in K[x]$.  Since an element of $K[x]$ lies
in $K[x^p]$ if and only if its derivative is zero, our hypotheses imply $g_i\in K[x^p]$.
Thus $g_i=\hat g_i(x)^p$ for some $\hat g_i\in\Kbar[x]$, so
$h=x^p\circ \hat g_1\circ f_1 = x^p\circ \hat g_2\circ f_2$, whence
$\hat g_1\circ f_1 = \hat g_2\circ f_2$.  In particular, $f_1$ and $f_2$ have a common
composite in $\Kbar[x]$ of degree less than $\deg(h)$, which contradicts
Theorem~\ref{algclosed}.
\end{proof}

\begin{remark}
Theorem~\ref{algclosed} was first proved by McConnell \cite{McConnell} in case $K$ is infinite,
and by Bremner and Morton \cite{BremnerMorton} in general.  Corollary~\ref{sep} is
due to Alexandru and Popescu \cite{AlexandruPopescu}.
\end{remark}

\section{Rational Composites and Polynomial Composites}
\label{sec rat}

\begin{theorem}
\label{ratpol}
If $f_1,f_2\in K[x]$ satisfy $K(f_1)\cap K(f_2)\ne K$, then $f_1$ and $f_2$ have a
common composite,
and moreover any minimal-degree common composite $h$ satisfies $K(f_1)\cap K(f_2)=K(h)$
and $K[f_1]\cap K[f_2]=K[h]$.
\end{theorem}

\begin{proof}
We use L\"uroth's theorem \cite[Thm.~2]{Schinzel}, which asserts that any subfield of $K(x)$
which properly contains $K$ must have the form $K(s)$.  Thus,
$K(f_1)\cap K(f_2)=K(\hat h)$ for some $\hat h\in K(x)$.
Write $\hat h=g_1\circ f_1=g_2\circ f_2$ with $g_1,g_2\in K(x)$, and
write $g_i=a_i/b_i$ with $a_i,b_i\in K[x]$ and $\gcd(a_i,b_i)=1$.
By inverting $\hat h,g_1,g_2$ if necessary, we may assume $\deg(a_1)\ge\deg(b_1)$.
Then
$$
        a_1(f_1(x))\cdot b_2(f_2(x)) = a_2(f_2(x))\cdot b_1(f_1(x)).
$$
In particular, $a_1(f_1(x))$ must divide the right hand side.  Since
$\gcd(a_1,b_1)=1$, some $K[x]$-linear combination of $a_1$ and $b_1$ equals 1;
substituting $f_1(x)$ for $x$ in this expression, it follows that 1 is a $K[x]$-linear
combination of $a_1\circ f_1$ and $b_1\circ f_1$, so $\gcd(a_1\circ f_1,b_1\circ f_1)=1$.
Thus, $a_1(f_1(x))$ divides $a_2(f_2(x))$.  By symmetry, they
must divide each other, so there is a constant $c$ such that 
$$
           a_1(f_1(x)) = c\cdot a_2(f_2(x)).
$$
In particular, $h_0:=a_1\circ f_1$ is in $K(f_1)\cap K(f_2)=K(\hat h)$.  But
$\deg(h_0)=\deg(\hat h)=[K(x)\col K(\hat h)]$, so in fact $K(h_0)=K(\hat h)$.

Now let $s$ be any common composite of $f_1$ and $f_2$.  Then $s\in K(f_1)\cap K(f_2)=K(h_0)$,
so $s=r\circ h_0$ with $r\in K(x)$.  It follows as above that $r\in K[x]$: write $r=a/b$ with
$a,b\in K[x]$ and $\gcd(a,b)=1$, so $\gcd(a\circ h,b\circ h)=1$, and since $b\circ h$ divides
$a\circ h$ we must have $\deg(b\circ h)=0$, whence $b$ is constant.  Thus
$K[f_1]\cap K[f_2]=K[h_0]$.  In particular, the minimal-degree common composites of
$f_1$ and $f_2$ are precisely the polynomials $\ell\circ h_0$ with $\ell\in K[x]$ of degree one.
The result follows.
\end{proof}

We can use this result to sharpen the conclusion of Theorem~\ref{algclosed}:

\begin{corollary}
\label{moncon}
If $h\in\Kbar[x]$ is a minimal-degree common composite of $f_1$ and $f_2$ over $\Kbar$,
then $\ell\circ h\in K[x]$ for some degree-one $\ell\in\Kbar[x]$.  In particular, if $h$ is monic
and has no constant term then $h\in K[x]$.
\end{corollary}

\begin{remark}
The anonymous referee informed us that, with some effort, one can prove
Corollary~\ref{moncon} via the linear algebra approach used to prove
Theorem~\ref{algclosed}.
\end{remark}

Another consequence of Theorem~\ref{ratpol} is that the study of common composites can be
reduced to the case where both polynomials have nonzero derivative:

\begin{corollary}
\label{insep}
For any $f_1,f_2\in K[x]$, write $f_i=\hat f_i\circ x^{p^{n_i}}$ with $n_i\ge 0$ and
$\hat f_i\in K[x]\setminus K[x^p]$, and suppose $n_1\ge n_2$.  
For any perfect field $L$ between $K$ and $\Kbar$, there exists $\tilde f_1\in L[x]\setminus
 L[x^p]$
such that $\hat f_1\circ x^{p^{n_1-n_2}} = x^{p^{n_1-n_2}}\circ\tilde f_1$.
Then $f_1$ and $f_2$ have a common composite over $K$
if and only if $\tilde f_1$ and $\hat f_2$ have a common composite over $L$.  Moreover, the
common composites of $f_1$ and $f_2$ over $L$ are precisely the polynomials of the form
$x^{p^{n_1-n_2}}\circ h\circ x^{p^{n_2}}$, where $h$ varies over the
common composites of $\tilde f_1$ and $\hat f_2$ over $L$.
\end{corollary}

\begin{proof}
  From the equation defining $\tilde f_1$, we see that $\tilde f_1$ is
  gotten from $\hat f_1$ by replacing each coefficient by its
  $p^{n_1-n_2}$-th root.  Thus $\tilde f_1\notin \Kbar[x^p]$.
  Next, since $f_1=x^{p^{n_1-n_2}}\circ \tilde f_1\circ x^{p^{n_2}}$
  and $f_2=\hat f_2\circ x^{p^{n_2}}$, the common composites of $f_1$
  and $f_2$ over $L$ are gotten by substituting $x^{p^{n_2}}$ into the
  common composites (over $L$) of $\hat f_2$ and $\bar
  f_1:=x^{p^{n_1-n_2}}\circ \tilde f_1$.  Write $q:=p^{n_1-n_2}$.  If
  $\tilde f_1$ and $\hat f_2$ have a common composite, then its
  $q$-th power is a composite of $\bar f_1$; thus $\bar f_1$ and
  $\hat f_2$ have a common composite if and only if $\tilde f_1$ and
  $\hat f_2$ do.  Hence $f_1$ and $f_2$ have a common composite over
  $L$ if and only if $\tilde f_1$ and $\hat f_2$ do, and by
  Theorem~\ref{algclosed} the former condition is equivalent to $f_1$
  and $f_2$ having a common composite over $K$.  So suppose $\tilde
  f_1$ and $\hat f_2$ have a common composite (over $L$), and let
  $\hat h$ be a common composite of minimal degree.  By
  Theorem~\ref{ratpol}, the common composites of $\tilde f_1$ and
  $\hat f_2$ are precisely the polynomials $\psi\circ \hat h$ with
  $\psi\in L[x]$.  Corollary~\ref{sep} implies that $\hat h\notin
  L[x^p]$.  Thus, $\psi\circ\hat h$ is in $L[x^q]$ if and only if
  $\psi\in L[x^q]$, or equivalently $\psi\circ\hat h\in L[\bar f_1]$.
  Hence the common composites of $\bar f_1$ and $\hat f_2$ are the
  polynomials $\varphi\circ x^q\circ\hat h$ with $\varphi\in L[x]$.
  Since $L$ is perfect, the set of $q$-th powers in $L[x]$
  equals $L[x^q]$, and the result follows.
\end{proof}

\begin{remark}
  The first two parts of Theorem~\ref{ratpol} were proved by Noether
  \cite{Noether} in the case of characteristic zero, and by McConnell
  \cite{McConnell} in general.  The third part of Theorem~\ref{ratpol}
  was proved by Schinzel \cite[Lemma 1, p.~18]{Schinzel}.
\end{remark}

We now give another proof of Theorem~\ref{ratpol} with a different flavor.

\begin{proof}[Second proof of Theorem~\ref{ratpol}]
First assume $K(x)$ is a separable extension of $K(f_1)\cap K(f_2)$.
By L\"uroth's theorem, $K(f_1)\cap K(f_2)=K(h)$ for some $h\in K(x)$.
By making a linear fractional change to $h$, we may assume that the infinite place of $K(x)$
lies over the infinite place of $K(h)$.
Let $N$ be the Galois closure of $K(x)/K(h)$, and let $G$, $H$, $A$, $B$ be
the subgroups of $\Gal(N/K(h))$ fixing $K(h)$, $K(x)$, $K(f_1)$, and $K(f_2)$.
Let $I$ be the inertia group in $N/K(h)$ of a place lying over the infinite place
of $K(h)$.  Then the corresponding inertia groups in $N/K(f_1)$ and $N/K(f_2)$
are $I\cap A$ and $I\cap B$.  Since $f_1$ and $f_2$ are polynomials, we have $A=H(I\cap A)$ and
$B=H(I\cap B)$.  For any subgroup $C$ of $G$, write $C_I$ for $I\cap C$.

Thus $H\langle A_I,B_I\rangle=\langle A_I,B_I\rangle H$, so $H\langle A_I,B_I\rangle$
is a group and thus equals $\langle A,B\rangle=G$.  Hence $HI=G$, so the infinite
place of $K(x)$ is the unique place of $K(x)$ lying over the infinite place of
$K(h)$, whence $h\in K[x]$.  Moreover, since the infinite place of $K(f_i)$ is
the unique place of $K(f_i)$ lying over the infinite place of $K(h)$, it follows that
$h$ is a common composite of $f_1$ and $f_2$.

Now assume $K(x)$ is an inseparable extension of $K(f_1)\cap K(f_2)$.
By L\"uroth's theorem, $\Kbar(f_1)\cap\Kbar(f_2)=\Kbar(h)$ for some $h\in\Kbar(x)$.
Write $f_1=x^{p^A}\circ\hat f_1$ and $f_2=x^{p^B}\circ\hat f_2$ where
$\hat f_i\in \Kbar[x]\setminus \Kbar[x^p]$.  Then there are
$g_i\in\Kbar(x)$ with $h=g_i\circ\hat f_i$.  Write
$g_1=x^{p^C}\circ \hat g_1$ and $g_2=x^{p^D}\circ \hat g_2$ with
$\hat g_i\in \Kbar(x)\setminus \Kbar(x^p)$.
Then $\Kbar(x)/\Kbar(\hat g_i\circ\hat f_i)$ is separable but
$\Kbar(\hat g_i\circ\hat f_i)/\Kbar(h)$ is purely inseparable,
so the latter extension is the maximal purely inseparable subextension
of $\Kbar(x)/\Kbar(h)$; in particular,
$\Kbar(\hat g_1\circ\hat f_1)=\Kbar(\hat g_2\circ\hat f_2)$.
Thus $\Kbar(x)$ is a separable extension of $\Kbar(\hat f_1)\cap\Kbar(\hat f_2)$,
so the result proved in the previous paragraphs implies that
$\Kbar(\hat f_1)\cap\Kbar(\hat f_2)=\Kbar(r)$ for some
$r\in\Kbar[x]$ which is a common composite of $\hat f_1$ and $\hat f_2$.
It follows easily that $\Kbar(f_1)\cap\Kbar(f_2)=\Kbar(r^{p^{\max(A,B)}})$.
Now Theorem~\ref{algclosed} implies that $f_1$ and $f_2$ have a common composite $\hat r\in K[x]$
with $\deg(\hat r)=\deg(r^{p^{\max(A,B)}})$, and since
$[K(x)\col K(f_1)\cap K(f_2)]\ge [\Kbar(x)\col \Kbar(f_1)\cap\Kbar(f_2)]$,
it follows that $K(f_1)\cap K(f_2)=K(\hat r)$ as desired.

We have shown that $K(f_1)\cap K(f_2)=K(h)$ where $h\in K[x]$ is a common composite of $f_1$
and $f_2$.  For any common composite $\hat h$ of $f_1$ and $f_2$, we have
$K(\hat h)\subseteq K(h)$, and moreover the infinite place of $K(h)$ is the unique place of $K(h)$
 lying
over the infinite place of $K(\hat h)$; thus $\hat h=r(h)$ for some $r\in K[x]$.
\end{proof}

This second proof generalizes at once to intersections of higher-genus function fields:

\begin{proposition}
Let $F$ be a finite extension of $K(x)$, and let $F_1$ and $F_2$ be subfields of $F$
which contain $K$.  Suppose $F$ is a finite separable extension of $F_0:=F_1\cap F_2$.
If a place $P$ of $F$ is totally ramified in both
$F/F_1$ and $F/F_2$, then $P$ is totally ramified in $F/F_0$.
\end{proposition}

In Section~\ref{sec comp} (following Theorem~\ref{comp})
we give a third proof of Theorem~\ref{ratpol}, which is a different type
of constructive proof.


\section{Degree Constraints}
\label{sec deg}

In this section we examine the possible degrees of common composites of $f_1$ and $f_2$.
By Theorem~\ref{ratpol}, the set of degrees of common composites equals the set of
multiples of some integer $n$, so it suffices to analyze $n$, which is the minimal degree
of any common composite.
Clearly any common composite has degree divisible by $\lcm(\deg(f_1),\deg(f_2))$.
Conversely, in characteristic zero we now show that if there is a common composite
then there is one of this minimal degree.  More generally this holds if there is a common
composite of degree not divisible by $p:=\charp(K)$:

\begin{theorem}
\label{lcm}
If $f_1,f_2\in K[x]$ have a common composite,
then they have a common composite of degree\/ $\lcm(\deg(f_1),\deg(f_2))p^s$ for some $s\ge 0$.
(Here we use the convention $0^0=1$.)
\end{theorem}

\begin{proof}
First assume $f_1,f_2\notin K[x^p]$.
Let $h(x)$ be a common composite of minimal degree.  Corollary~\ref{sep} and Theorem~\ref{ratpol}
imply that $K(x)/K(h(x))$ is separable and $K(h)=K(f_1)\cap K(f_2)$.
Let $L$ be the Galois closure of $K(x)/K(h(x))$, and let $G,A,B,H$ be the subgroups of
$\Gal(L/K(h(x)))$ fixing $h(x)$, $f_1(x)$, $f_2(x)$, and $x$.  Then $G=\langle A,B\rangle$.
Let $P$ be a place of $L$ lying over the infinite place of $K(h(x))$,
and let $I$ be the inertia group of $P$ in $L/K(h(x))$.  Since $h$ is a polynomial,
$G=HI$.

For any group $C$ with $H\le C\le G$, let $C_I:=C\cap I$.  Clearly $C$
contains $HC_I$, and since $G=HI$ we have $C=HC_I$.  Moreover,
$[C\col H]=[C_I\col H_I]$.  Since $H\langle A_I,B_I\rangle=\langle
A_I,B_I\rangle H$, the set $H\langle A_I,B_I\rangle$ is a group and
thus equals $\langle A,B\rangle=G$.  Hence $I=G_I=\langle
A_I,B_I\rangle$.  Recall the structure of inertia groups (cf., e.g.,
\cite[Cor.~4 to Prop.~7, \S~IV.2]{Serre}): $I$ is the
semidirect product $V\rtimes D$ where $V$ is a normal $p$-subgroup
and $D$ is cyclic of order not divisible by~$p$.  Since $I=\langle
VA_I,VB_I\rangle$, we have $I/VH_I=\langle VA_I/VH_I,VB_I/VH_I\rangle$, and
since these are cyclic groups we see that $[I\col VH_I]$ is the least
common multiple of $[VA_I\col VH_I]$ and $[VB_I\col VH_I]$.  Finally,
$\deg(h)=[I\col H_I]$ and $\deg(f_1)=[A_I\col H_I]$, so $[V\col VH_I]$ and
$[VA_I\col VH_I]$ are the maximal divisors of $\deg(h)$ and $\deg(f_1)$
which are not divisible by~$p$, whence
$\deg(h)=\lcm(\deg(f_1),\deg(f_2))p^t$ with $t\ge 0$.

Now for arbitrary $f_1,f_2\in K[x]$ having a common composite,
Corollary~\ref{insep} implies that the minimal degree of any common
composite of $f_1$ and $f_2$ over $\Kbar$ is a power of $p$ times the
minimal degree of any common composite of two related polynomials
$\tilde f_1,\hat f_2\in \Kbar[x]\setminus\Kbar[x^p]$, where both
$\deg(f_1)/\deg(\tilde f_1)$ and $\deg(f_2)/\deg(\hat f_2)$ are powers
of~$p$.  Since $\tilde f_1,\hat f_2\notin \Kbar[x^p]$, it follows from
above that the minimal degree of any common composite of $\tilde f_1$
and $\hat f_2$ is $\lcm(\deg(\tilde f_1),\deg(\hat f_2))p^s$ for some
$s\ge 0$, so the minimal degree of any common composite of $f_1$ and
$f_2$ over $\Kbar$ is $\lcm(\deg(f_1),\deg(f_2))p^t$ with $t\ge 0$.
The result now follows from Theorem~\ref{algclosed}.
\end{proof}

\begin{remark}
  Theorem~\ref{lcm} was proved by Engstrom \cite{Engstrom} in the case of
  characteristic zero, and his proof extends at once to the case
  where $f_1$ and $f_2$ have a common composite of degree not
  divisible by $p$ (cf.\ \cite[Thm.~5]{Schinzel}).  This elegant
  proof is completely different from ours (for instance it depends on
  nothing beyond the division algorithm in $K[x]$), and it would be
  interesting to try to extend Engstrom's argument to prove our full
  result.  In case $f_1$ and $f_2$ have a common composite of degree
  not divisible by $p$, our proof is essentially a modernized account
  of an argument due to Ritt \cite{Ritt}; an alternate treatment of
  Ritt's proof in this case, using fields and power series instead of
  groups and inertia groups, is in \cite{McConnell}.  The basic ideas
  in \cite{McConnell} can be discerned by scrutinizing the proof of
  \cite[Thm.~3.6]{FriedMacRae}, though significant effort is
  required since the latter proof contains errors in nearly every
  line.  An incorrect generalization of Theorem~\ref{lcm} is given as
  \cite[Thm.~2.1]{AlexandruPopescu}; specifically, they assert that
  the result for degrees also holds for the ramification indices under
  any prescribed place of $K(x)$.  A counterexample is $f_1=x^2$ and
  $f_2=x^3-x$ over $K=\bC$ at the place $x=1$, since $x=1$ is
  unramified in $K(x)/K(f_1)$ and $K(x)/K(f_2)$ but ramifies in
  $K(x)/(K(f_1)\cap K(f_2))=K(x)/K((x^3-x)^2)$.  The mistake in the
  proof of \cite[Thm.~2.1]{AlexandruPopescu} is the assertion that
  the completions of $K(f_1)$ and $K(f_2)$ (at places under the
  prescribed place of $K(x)$) intersect in the completion of
  $K(f_1)\cap K(f_2)$, which is not generally true.
\end{remark}

It is not possible to remove the power of $p$ from the conclusion of Theorem~\ref{lcm}.
For example, there are polynomials over $\F_2$ of degrees 11 and 13 whose
least-degree common composite has degree $143\cdot 2^{60}$, and there are
polynomials over $\F_2$ of degrees 1447 and 1451 whose least-degree common 
composite has degree $1447\cdot 1451\cdot 2^{1048350}$.
These are special cases of the following result.

\begin{proposition}
\label{funny}
Suppose $p:=\charp(K)$ is nonzero.  If $p\nmid n$ then a minimal-degree
common composite of $x^n$ and $x^{p^r}-x$ is $(x^{p^{rd}}-x)^n$, where $d$ is the
multiplicative order of $p^r \bmod n$.  If $p\nmid nm$ and $n,m>1$ then a
minimal-degree
common composite of $x^n$ and $(x-1)^m$ is $(x^{p^d}-x)^{\lcm(m,n)}$, where
$d$ is the multiplicative order of $p \bmod \lcm(m,n)$.
\end{proposition}

\begin{proof}
Let $\Kbar$ be an algebraic closure of $K$, and let $\zeta$ be a primitive
$n$-th root of unity in $\Kbar$.  Then $\Kbar(x)/\Kbar(x^n)$ is Galois with group
generated by $\sigma:x\mapsto \zeta x$, and $\Kbar(x)/\Kbar(x^{p^r}-x)$ is Galois with group
$H$ consisting of the various maps $x\mapsto x+\alpha$ with $\alpha\in\F_{p^r}$.
The subgroup $G$ of $\Aut_{\Kbar} \Kbar(x)$ generated by $\sigma$ and $H$ consists of
the maps $x\mapsto \mu x + \nu$ where $\mu\in\langle\zeta\rangle$ and $\nu\in\F_{p^r}(\zeta)$.
Here $\#G=np^{rd}$, where $d:=[\F_{p^r}(\zeta)\col \F_{p^r}]$ is the multiplicative order
of $p^r$ mod $n$.  The group $G$ fixes $h(x):=(x^{p^{rd}}-x)^n$, so since
$\deg h=\#G$ we see that $\Kbar(h)$ is the subfield of $\Kbar(x)$ fixed by $G$,
whence $\Kbar(h)$ is the intersection of $\Kbar(x^n)$ and $\Kbar(x^{p^r}-x)$.  This shows
that $h$ is a minimal-degree common composite of $x^n$ and $x^{p^r}-x$ over $\Kbar$.
Since $h\in K[x]$, it is also a minimal-degree common composite over $K$.

Now let $\eta$ be a primitive $m$-th root of unity in $\Kbar$.  Then the extension
$\Kbar(x)/\Kbar((x-1)^m)$ is Galois with group generated by $\gamma:x\mapsto 1+\eta(x-1)$.
Let $H$ be the subgroup of $\Aut_{\Kbar} \Kbar(x)$ generated by $\sigma$ and $\gamma$.
Then $H$ contains the commutator $\gamma^{-1}\sigma^{-1}\gamma\sigma:x\mapsto
x-(\eta-1)(\zeta-1)$.
One easily checks that $H$ consists of the maps $x\mapsto\mu x+\nu$ where
$\mu\in\langle\zeta,\eta\rangle$ and $\nu\in\F_p(\zeta,\eta)$.
Moreover, $H$ fixes $j(x):=(x^{p^d}-x)^{\lcm(m,n)}$, where $d$ is the multiplicative
order of $p$ mod $\lcm(m,n)$.  Since $\deg(j)=\#H$, it follows as above that
$j$ is a minimal-degree common composite of $x^n$ and $(x-1)^m$ over $\Kbar$,
and hence over $K$.
\end{proof}

\begin{remark}
The second part of Proposition~\ref{funny} was first proved by Bremner and
 Morton \cite{BremnerMorton}.
\end{remark}

In case $\charp(K)=0$, Theorem~\ref{lcm} says that if $f_1,f_2\in
K[x]$ have a common composite then they have one of degree
$\lcm(\deg(f_1),\deg(f_2))$.  Proposition~\ref{funny} shows that this
is no longer true when $\charp(K)>0$.  Specifically, for any prime
$p$, the least-degree common composite of $x^2$ and $x^2-x$ over
$\F_p$ is $(x^p-x)^2$.  Thus, two degree-$2$ polynomials can have
lowest-degree common composite of arbitrarily large degree.  We make
the following definition to give a framework for recovering some
analogue of the characteristic zero result, by restricting to
polynomials over a fixed field, or polynomials over fields of a fixed
characteristic.

\begin{definition}
Given integers $n_1,n_2>1$ and a field $K$, let $N(n_1,n_2,K)$ be the supremum of the integers
$r(f_1,f_2,K)$, where
\begin{itemize}
\item $f_1,f_2\in K[x]$ have a common composite and satisfy $\deg(f_1)=n_1$ and $\deg(f_2)=n_2$,
and
\item $r(f_1,f_2,K)$ is the lowest degree of any common composite of $f_1$ and $f_2$.
\end{itemize}
For any prime number $p$ (and for $p=0$), let $N(n_1,n_2,p)$ be the supremum of the values
$N(n_1,n_2,K)$,
where $K$ varies over all fields of characteristic~$p$.
\end{definition}

Theorem~\ref{lcm} implies that $N(n_1,n_2,0)=\lcm(n_1,n_2)$, and more generally that
$r(f_1,f_2,K)=\lcm(n_1,n_2)\charp(K)^s$.
However, Proposition~\ref{funny} shows that in positive characteristic there are
examples with arbitrarily large $s$.  But our examples have $n_1+n_2\to\infty$, and we do
not know whether one can bound $s$ in terms of $n_1$, $n_2$ and $K$, or even just in terms
of $n_1$ and $n_2$.  In fact, every example we know (when $p=\charp(K)>0$) satisfies
$s\le\lcm(n_1,n_2)$.  We now prove that, when $n_1=n_2=2$, we can actually take $s\le 1$:

\begin{proposition}
\label{deg2}
Any two degree-$2$ polynomials over a field of characteristic~$p>0$ have a common composite
of degree $2p$.
\end{proposition}

\begin{proof}
Let $K$ be a field of characteristic~$p$, and let $f_1$ and $f_2$ be degree-$2$ polynomials
in $K[x]$.
If $K(x)/K(f_1)$ is not separable, then $p=2$ and $f_1=ax^2+b$, so $f_2^2$ is a common composite
of $f_1$ and $f_2$ of degree $2p$.  Henceforth assume $K(x)/K(f_1)$ and $K(x)/K(f_2)$ are
separable.  Thus these extensions
are Galois.  Moreover, writing $f_1=ax^2-bx+c$, the Galois group of $K(x)/K(f_1)$ is generated by
$x\mapsto b/a-x$.  Thus, $\Gal(K(x)/K(f_1))$ and $\Gal(K(x)/K(f_2))$ are generated by
$\sigma_1:x\mapsto \alpha_1-x$ and
$\sigma_2:x\mapsto \alpha_2-x$, for some $\alpha_1,\alpha_2\in K$.
Now, $K(f_1)\cap K(f_2)$ is the subfield of $K(x)$ fixed by $H:=\langle \sigma_1,
\sigma_2\rangle$.
Since $\sigma_1$ and $\sigma_2$ have order 2, they generate a dihedral group of order
 twice
the order of the composite map $\sigma_1\sigma_2:x\mapsto (\alpha_1-x)\circ (\alpha_2-x)=
\alpha_1-\alpha_2+x$.
Since the latter map has order $1$ or $p$, it follows that $\#H\mid 2p$.  Now the result follows
from Theorem~\ref{ratpol}.
\end{proof}

The anonymous referee suggested the following alternate proof:
\begin{proof}[Second proof of Proposition~\ref{deg2}]
Let $K$ be a field of characteristic $p$, and let $f_1$ and $f_2$ be degree-$2$
polynomials in $K[x]$.  By replacing $f_i$ with $\ell_i\circ f_i$ for a suitable
degree-$1$ polyomial $\ell_i\in K[x]$, we may assume $f_1=x^2+ax$ and $f_2=x^2+bx$.
If $a=b$ then $f_1$ is already a common composite; hence we assume $a\ne b$.
Then $f_1$ and $f_2$ have a common composite of degree at most $2n$ if and only if
the polynomials $1,f_1,f_2,f_1^2,f_2^2,\dots,f_1^n,f_2^n$ are linearly dependent.  These
polynomials span the same space as the polynomials
$1,f_1-f_2,f_2,f_1^2-f_2^2,f_2^2,\dots,f_1^n-f_2^n,f_2^n$.
Since the leading term of
$f_2^i$ is $x^{2i}$, and the leading term of $f_1^i-f_2^i$ is $i(a-b)x^{2i-1}$,
the matrix of coefficients of these polynomials is triangular, and it has no zero
entries on the main diagonal if and only if $n<p$.  Thus $f_1$ and $f_2$ have a
common composite of degree $2p$, and no common composite of lower degree.
\end{proof}

It would be interesting to determine further values of $N(n_1,n_2,p)$, or even to determine
whether these values are finite.
One can attempt to produce infinite values of $N(n_1,n_2,p)$ by modifying the proof of
 Theorem~\ref{lcm}.
Below is a group-theoretic example
satisfying the conditions used in that proof, such that $H_I,A_I,B_I,I$ have orders
$1,2,3,2\cdot 3^{2n+1}$ respectively, where $n$ can be any positive integer.  If this
group-theoretic setup could be realized by polynomials $f_1$ and $f_2$ in
characteristic~3, then there would be polynomials of degrees~2 and~3
whose lowest-degree common composite has degree $2\cdot 3^{2n+1}$.

\begin{exampleme}
Let $I$ be the group generated by $a,b,c$ subject to the relations
$b^{3^n}=c^{3^n}=a^6=1$, $bc=cb$, $a^{-1}ba=c^{-1}$, $a^{-1}ca=bc$.
One can check that $I=VC$ where $V=\langle b,c,a^2\rangle$ has order $3^{2n+1}$
and $C=\langle a^3\rangle$ has order $2$.
Now $\langle a^3b,a^2b\rangle$ contains $a$ and $b$ and
thus contains $c=a^{-1}b^{-1}a$.  Hence $I=\langle A_I,B_I\rangle$ where
$A_I:=\langle a^3b\rangle$ and $B_I:=\langle a^2b\rangle$.  Finally, one can check
that $\#A_I=2$ and $\#B_I=3$.
\end{exampleme}

\noindent
In a subsequent paper we will show that the above configuration does not happen,
and in fact we will compute $N(2,3,p)$.  However, our proof uses a different
framework, and does not give a simple explanation why the above configuration doesn't
occur.  It would be interesting to know a general constraint on the Galois groups
associated to a polynomial which would preclude this setup from being realizable.

Theorem~\ref{lcm} is a `least common multiple' result.  We now give a companion
`greatest common divisor' result.

\begin{theorem}
\label{gcd}
If $f_1,f_2\in K[x]$ have a common composite of degree not divisible by $p$,
then there are $g_1,g_2,r\in K[x]$ with $\deg(r)=\gcd(\deg(f_1),\deg(f_2))$
such that $f_1=g_1\circ r$ and $f_2=g_2\circ r$.
\end{theorem}

\begin{proof}
We use the notation from the proof of Theorem~\ref{lcm}.  Thus $H(A_I\cap B_I)=
(A_I\cap B_I)H$, so $H(A_I\cap B_I)$ is a group, and equals $A\cap B$.
By L\"uroth's theorem, the subfield of $L$ fixed by $H(A_I\cap B_I)$ has the
form $K(r(x))$ for some rational function $r(x)$.  By making a linear fractional
change to $r(x)$ if necessary, we may assume that the infinite place of $K(r(x))$
lies under the infinite place of $K(x)$.  Since the latter place is totally
ramified in $K(x)/K(r(x))$, it follows that $r(x)$ is a polynomial.  Moreover,
the infinite place of $K(r(x))$ is the unique place lying over the infinite place
of $K(f(x))$, so $f_1=g_1\circ r$ for some polynomial $g_1$, and likewise
$f_2=g_2\circ r$.

It remains only to determine the degree of $r$, which equals
$[H(A_I\cap B_I)\col H]=\#(A_I\cap B_I)$.
Since $A_I$ and $B_I$ are subgroups of the cyclic group $I$, we have
$\#(A_I\cap B_I)=\gcd(\#A_I,\#B_I)$.  Thus the degree of $r$ is
$\gcd([A\col H],[B\col H])=\gcd(\deg(f_1),\deg(f_2))$.
\end{proof}

\begin{remark}
Theorem~\ref{gcd} was proved by Engstrom \cite{Engstrom} in the case of
characteristic zero, and his proof extends at once to the general case
(cf.\ \cite[Thm.~5]{Schinzel}).  The situation is the same as for
Theorem~\ref{lcm}: Engstrom's argument (as simplified by Schinzel)
uses just polynomials, and no Galois theory.  Our Galois-theoretic proof is
a modernized version of an argument of Ritt's \cite{Ritt}, and a complicated
field-theoretic version of Ritt's argument (with numerous errors) is in \cite{FriedMacRae}.
\end{remark}

Note that the hypothesis on the degree in Theorem~\ref{gcd} is
necessary---for instance, if $f_1=x^2+ax$ and $f_2=x^2+bx$ with $a\ne b$, then certainly there
is no $r$ satisfying the conclusion of Theorem~\ref{gcd}, but Proposition~\ref{deg2}
says that $f_1$ and $f_2$ have a common composite over any field of positive characteristic.

Theorems \ref{lcm} and \ref{gcd} show that the existence of a common composite
of degree not divisible by $p$ is a very unusual occurrence.  For instance,
if polynomials $f_1$ and $f_2$ of the same degree have such a common composite, then
$f_1=\ell\circ f_2$ for some degree-$1$ polynomial $\ell$.


\section{The Tame Case}
\label{sec tame}

In this section we describe all pairs of polynomials $f_1,f_2\in K[x]$ which
have a common composite of degree not divisible by $\charp(K)$.  The statement of
the result involves the Dickson polynomials, which are defined as follows.  For any
$\alpha\in K$ and $n>0$, define $D_n(x,\alpha)\in K[x]$
by
$$
D_n(x,\alpha)=\sum_{i=0}^{\lfloor n/2\rfloor} \frac{n}{n-i} \binom{n-i}i (-\alpha)^i x^{n-2i}.
$$
The key property of Dickson polynomials is that $D_n(x+\alpha/x,\alpha)=x^n+(\alpha/x)^n$.

\begin{theorem}
\label{char0}
Suppose $f_1,f_2\in K[x]$ satisfy $\deg(f_1)\ge\deg(f_2)>1$ and $\charp(K)\nmid\deg(f_1)\deg(f_2)$.
Then $f_1$ and $f_2$ have a common composite of degree not divisible by $\charp(K)$ if and only if
there are degree-$1$ polynomials $\ell_1,\ell_2\in K[x]$ and a polynomial
$h(x)\in K[x]$ of degree $\gcd(\deg(f_1),\deg(f_2))$ such that either
\begin{enumerate}
\item[\rm 1.] $f_1=\ell_1\circ x^r P(x^n) \circ h(x)$ and $f_2=\ell_2\circ x^n\circ h(x)$,
   where $r,n>0$ and $P\in K[x]$; or
\item[\rm 2.] $f_1=\ell_1\circ D_m(x,\alpha)\circ h(x)$ and $f_2=\ell_2\circ D_n(x,\alpha)\circ
   h(x)$, where $\alpha\in K$ and $m,n>0$.
\end{enumerate}
\end{theorem}

\begin{proof}
Suppose $f_1$ and $f_2$ have a common composite of degree not divisible by $p:=\charp(K)$.
By Theorem~\ref{gcd},
there are $g_1,g_2,h\in K[x]$ such that $f_i=g_i\circ h$ and $\deg(h)=
\gcd(\deg(f_1),\deg(f_2))$.
By Theorem~\ref{lcm}, $g_1$ and $g_2$ have a common composite of degree $\lcm(\deg(g_1),
\deg(g_2))$.
Now the result follows from Theorem~\ref{ritt2} below.
\end{proof}

In the following result, if $\ell$ is a degree-$1$ polynomial over a field $K$,
we write $\ell^{\langle -1\rangle}$ to denote the functional inverse of $\ell$;
thus $\ell^{\langle -1\rangle}$ is the unique degree-$1$ polynomial over $K$ for
which $\ell^{\langle -1\rangle}(\ell(x))=x$, or equivalently
$\ell(\ell^{\langle -1\rangle}(x))=x$.

\begin{theorem}[Zannier]
\label{ritt2}
Suppose $a,b,c,d\in K[x]$ satisfy $\deg(a)=\deg(d)=m>1$ and $\deg(b)=\deg(c)=n>1$,
where $\gcd(m,n)=1$ and $m>n$ and $a'c'\ne 0$.  Then $a(b)=c(d)$ holds if and only if
there are degree-$1$ polynomials $\ell_1,\ell_2,\ell_3,\ell_4\in K[x]$ such that either
\begin{enumerate}
\item[\rm 1.] $\ell_1\circ a\circ \ell_3^{\langle -1\rangle}= x^r P(x)^n$ and $\ell_3\circ b\circ \ell_2=x^n$
  and $\ell_1\circ c\circ \ell_4^{\langle -1\rangle}= x^n$ and $\ell_4\circ d\circ \ell_2=x^r P(x^n)$,
  where $P\in K[x]$ and $r=m-n\deg(P)>0$; or
\item[\rm 2.] $\ell_1\circ a\circ \ell_3^{\langle -1\rangle}= D_m(x,\alpha^n)$ and $\ell_3\circ b\circ \ell_2=
 D_n(x,\alpha)$
  and $\ell_1\circ c\circ \ell_4^{\langle -1\rangle}= D_n(x,\alpha^m)$ and $\ell_4\circ d\circ \ell_2=
  D_m(x,\alpha)$, where $\alpha\in K$.
\end{enumerate}
\end{theorem}

\begin{remark}
Theorem~\ref{ritt2} was proved by Zannier \cite{Zannier}; an alternate exposition of his
proof is in~\cite[Thm.~8]{Schinzel}.  Previously special cases had been proved by
Ritt~\cite{Ritt},
Levi~\cite{Levi}, Dorey and Whaples~\cite{DoreyWhaples}, Schinzel~\cite{Schinzel2}, and
Tortrat~\cite{Tortrat}.
\end{remark}

Theorem~\ref{char0} shows in a strong sense that, when $\charp(K)=0$, very few pairs of
 polynomials $(f_1,f_2)$
have a common composite.  We suspect that the same qualitative behavior occurs in positive
characteristic, but it is difficult to prove significant results in this direction.


\section{Fiber-finding}
\label{sec fib}

In this section we give two algorithms which produce either a common composite or
a proof that there is no such of degree less than a prescribed bound.
The idea of the first algorithm is simple: if $f_1$ and $f_2$ have a common composite $h$,
then any $\alpha,\beta\in\Kbar$ with $f_1(\alpha)=f_1(\beta)$ also satisfy
$h(\alpha)=h(\beta)$.  Thus, starting with some $\alpha\in \Kbar$, we compute
all $\beta\in\Kbar$ with $f_1(\alpha)=f_1(\beta)$.  Then for each $\beta$
we compute all $\gamma\in\Kbar$ with $f_2(\beta)=f_2(\gamma)$.  Note that
$h(\gamma)=h(\beta)=h(\alpha)$.  Continuing this process, we find more and more
elements of $\Kbar$ which have the same $h$-value.  This gives a lower bound
on the degree of $h$; conversely, we show in Section \ref{sec comp} that,
if this process produces only finitely many elements of $\Kbar$, then we can
determine whether $f_1$ and $f_2$ have a common composite.

In the second algorithm we work with polynomials over $K$ rather than elements of $\Kbar$.
In this case it is convenient to assume the $f_i$ are monic.
Suppose we have a nonconstant $r\in K[x]$ which divides
our hypothesized (minimal degree) common composite $h$.  Let $m$ be the minimal polynomial
of $f_i$ mod $r$, i.e., $m$ is the minimal degree monic polynomial in $K[x]$
such that $m\circ f_i$ is divisible by $r$.  Then $m\circ f_i$ divides $h$.
By iterating this process, we can quickly build up large-degree factors of $h$.
We can start this process with $r_0=x$.  After one step, we have $r_1=(x-f_1(0))\circ f_1
= f_1 - f_1(0)$.  The polynomials $r_j$ alternate between composites of
$f_1$ and composites of $f_2$.  Therefore, if this process ever
stabilizes (by giving $r_j=r_{j+1}$ for some $j>0$), then
the final $r_j$ is a minimal degree common composite of $f_1$ and $f_2$.

\begin{exampleme}
\label{char3}
Let $f_1=x^2$ and $f_2=x^3+x^2+x$, where $\charp(K)=3$.  Then we start with
$r_1:=f_1=x^2$.
The minimal polynomial of $f_2$ mod $r_1$ is $x^2$, so we put
$r_2:=x^2\circ f_2=x^6-x^5-x^3+x^2$.  The minimal polynomial of $f_1$ mod $r_2$ is
$x^5-x^4-x^2+x$,
so we put $r_3:=(x^5-x^4-x^2+x)\circ f_1=x^{10}-x^8-x^4+x^2$.
The minimal polynomial of $f_2$ mod $r_3$ is $m:=x^6+x^5+x^3+x^2$,
so we put $r_4:=m\circ f_2=x^{18}-x^{14}-x^6+x^2$.
Finally, the minimal polynomial of $f_1$ mod $r_4$ is $x^9-x^7-x^3+x$, and
$r_4=(x^9-x^7-x^3+x)\circ f_1$, so $r_4$ is a minimal-degree common composite.
\end{exampleme}

The above algorithms are actually two incarnations of the same idea.
In the first algorithm we explore the fiber $\{\zeta:h(\zeta)=h(\alpha)\}$.
Letting $Z$ be the set of $\zeta$'s seen up to a given step, we can put
$r:=\prod_{\zeta\in Z}(x-\zeta)$.  We know $r(x)$ divides $h(x)-h(\alpha)$;
we may assume $h(\alpha)=0$, so $r$ divides $h$.  Suppose the next step involves
equating values of $f_i$, and let $\hat Z$ be the next set of $\zeta$'s.
Let $v(x)$ be obtained by eliminating $z$ from the system:
\begin{eqnarray*}
r(z)&=&0\\
f_i(z)&=&f_i(x)
\end{eqnarray*}
(i.e., $v(x)$ generates the intersection of the ideal $(r(z),f_i(z)-f_i(x))$ with $K[x]$).
Then every root of $v$ lies in $\hat Z$, and
every element of $\hat Z$ is a root of $v$.  It turns out that
$v=m\circ f_i$, where $m$ is the minimal polynomial of $f_i$ mod $r$.
Thus, both our algorithms produce the same set of elements of $\Kbar$ at each step;
the main difference between them is that the second algorithm keeps track of multiplicities,
while the first does not.

Here is an example where the second algorithm can be used to prove that two polynomials
have no common composite.

\begin{exampleme}\label{strange}
Let $f_1=x^2-x$ and $f_2=x^3-x^2$.  We start with $r_1:=f_1$.
Inductively, we show that $r_{2j+1}=f_1^{2^j}$ and $r_{2j+2}=f_2^{2^j}$.
Indeed, if $r_{2j+1}=f_1^{2^j}$ then its roots $x=0$ and $x=1$ each have
multiplicity $2^j$; since $x=0$ and $x=1$ are roots of $f_2$ of multiplicities
$2$ and $1$, it follows that $r_{2j+2}=f_2^{2^j}$.  Thus the roots of $r_{2j+2}$
are again $x=0$ and $x=1$, this time with multiplicities $2^{j+1}$ and $2^j$;
since $x=0$ and $x=1$ are simple roots of $f_1$, it follows that
$r_{2j+3}=f_1^{2^{j+1}}$.  Since the degrees of the $r_j$ grow without bound,
$f_1$ and $f_2$ have no common composite.
\end{exampleme}

If we apply the first algorithm with $\alpha=0$ to the polynomials in the above example,
we quickly find a stable set $Z=\{0,1\}$.
This example is better understood in the context of the next two sections: the first
algorithm terminates with $Z=\{0,1\}$ because that set is {\em
compatible}\/ (see Section~\ref{sec comp}).  The second algorithm fails to terminate
because the set $Z$ is {\em inconsistent}\/ (see Example~\ref{firstexample}).

While we suspect that the above example illustrates a rare situation,
it is worth modifying the second algorithm so that, if the set $Z$ of
roots of $r$ stabilizes, we check $Z$ for consistency.

\section{Compatible Consistent Sets}
\label{sec comp}

Let $f_1$ and $f_2$ be nonconstant polynomials over $K$.  If $f_1$ and $f_2$ have a
common composite $h$ then, for any $\alpha\in \Kbar$, the $h$-fiber
$\{\beta\in\Kbar: h(\beta)=h(\alpha)\}$ is
a finite subset of $\Kbar$ which is simultaneously a union of $f_1$-fibers and
a union of $f_2$-fibers.
We generalize this to arbitrary $f_1$ and $f_2$ (which might not have a common
composite) as follows:

\begin{definitionme}
A nonempty finite subset of $\Kbar$ is \emph{compatible} if it is simultaneously
a union of $f_1$-fibers and a union of $f_2$-fibers.
\end{definitionme}

We will show that, if there is a compatible set, then there is a common composite
precisely when a certain easily checkable condition is met.  To motivate this
extra condition, assume again that $f_1$ and $f_2$ have a common composite $h$.
For each $a\in \Kbar$, let $\ell(a)$ be the ramification
index of $x=a$ in the extension $\Kbar(x)/\Kbar(h(x))$; in other words, $\ell(a)$ is the
multiplicity of $x=a$ as a root of $h(x)-h(a)$.  Likewise, let $m_i(a)$
be the ramification index of $x=a$ in the extension $\Kbar(x)/\Kbar(f_i(x))$.  Then
$m_i(a)$ divides $\ell(a)$, and moreover if $a,b\in \Kbar$ satisfy $f_i(a)=f_i(b)$ for
some $i$ then the ramification index of $f_i(x)=f_i(a)$ in $\Kbar(f_i(x))/\Kbar(h(x))$ is
\begin{equation*}
\frac{\ell(a)}{m_i(a)}=\frac{\ell(b)}{m_i(b)}.
\end{equation*}
In general, when $f_1$ and $f_2$ are not assumed to have a common composite,
we make the following definition.  Again, $m_i(a)$ is the ramification index of $x=a$
in the extension $\Kbar(x)/\Kbar(f_i(x))$.

\begin{definitionme}
\label{cons}
A subset $A\subseteq\Kbar$ is \emph{consistent} if there is a function $\ell$ on $A$
such that
\begin{enumerate}
\item for each $a\in A$, $\ell(a)$ is a positive integer multiple of both $m_1(a)$
and $m_2(a)$; and
\item for $a,b\in A$ and $i\in\{1,2\}$, if $f_i(a)=f_i(b)$ then
$\ell(a)/m_i(a)=\ell(b)/m_i(b)$.
\end{enumerate}
\end{definitionme}

The above discussion implies
\begin{proposition}
\label{conscomp}
If $f_1$ and $f_2$ have a common composite of degree $n$, then every element of $\Kbar$ is
contained in a compatible set of size at most~$n$, and every subset of $\Kbar$ is consistent
via the labeling defined by the ramification index in $\Kbar(x)/(\Kbar(f_1)\cap\Kbar(f_2))$.
\end{proposition}

We now prove a converse result, which implies Theorem~\ref{thmintro}:

\begin{theorem}
\label{comp}
If there is a compatible consistent set $A$, then $f_1$ and $f_2$
have a common composite over $K$.  Explicitly, if $\ell:A\to\Z$ is a consistent labeling on $A$, then
$h:=\prod_{a\in A}(x-a)^{\ell(a)}$ is a common composite over $\Kbar$.
\end{theorem}

\begin{proof}
We may assume $f_1$ and $f_2$ are monic.  Let $B_1=\{f_1(a):a \in A\}$, and for
each $b\in B_1$ pick an element $a_b\in A$ with $f_1(a_b)=b$.  Let $A_1=\{a_b:b \in B_1\}$.
Now we compute
\begin{align*}
\prod_{a\in A}(x-a)^{\ell(a)} &= \prod_{\hat a\in A_1}\prod_{\substack{a \in A\\f_1(a)=
f_1(\hat a)}} (x-a)^{\ell(a)} \\
&=\prod_{\hat a\in A_1}\Bigl(\prod_{\substack{a \in A\\f_1(a)=f_1(\hat a)}}(x-a)^{m_1(a)}\Bigr)^{
\ell(\hat a)/m_1(\hat a)} \\
&=\prod_{\hat a\in A_1}\left(f_1(x)-f_1(\hat a)\right)^{\ell(\hat a)/m_1(\hat a)}\\
&=\Bigl(\prod_{\hat a\in A_1}\left(x-f_1(\hat a)\right)^{\ell(\hat a)/m_1(\hat a)}\Bigr) \circ
 f_1(x),
\end{align*}
where the two middle equalities hold because $A$ is consistent and compatible, respectively.
Thus, the polynomial $h:=\prod_{a\in A}(x-a)^{\ell(a)}$ is a composite of $f_1$ over $\Kbar$; but
 likewise it
is a composite of $f_2$ over $\Kbar$, so it is a common composite over $\Kbar$.
It follows by Theorem~\ref{algclosed} that $f_1$ and $f_2$ have a common composite over $K$.
\end{proof}

This result has several consequences.  For one thing, it gives yet another proof of the first
part of Theorem~\ref{ratpol}, namely that $K(f_1)\cap K(f_2)\ne K$ implies $f_1$ and $f_2$ have a
common composite: for in this case $K(f_1)\cap K(f_2)=K(h)\ne K$, and the proof of
Proposition~\ref{conscomp} shows there are compatible consistent subsets of $\Kbar$.
More importantly, in Theorem~\ref{comp} we exhibited a specific common composite $h$ over $\Kbar$.
The shape of this polynomial $h$ enables us to control the ramification in a minimal-degree common
composite in terms of the ramification in $f_1$ and $f_2$; in a subsequent paper we will show how
this can be used to prove that two polynomials have no common composite.

\begin{corollary}
\label{ramcor}
If $f_1,f_2\in K[x]$ have a common composite, then the ramification index of $x=a$ in
$\Kbar(x)/(\Kbar(f_1)\cap \Kbar(f_2))$ is a divisor of $\ell(a)$, for any consistent labeling
$\ell$ on any compatible set containing $a$.
\end{corollary}

Another consequence of Theorem~\ref{comp} is a description of the minimal compatible sets,
in case there is a common composite.  We need a lemma before stating the result:
\begin{lemma}
If $f_1$ and $f_2$ have a common composite, and $A\subset\Kbar$ is a minimal compatible set,
then there is a consistent labeling $\ell_0:A\to\Z$ such that every consistent labeling
$\ell:A\to\Z$ has the form $\ell=n\ell_0$ with $n$ a positive integer.
\end{lemma}
\begin{proof}
Pick some $a\in A$ and some consistent labeling $\ell:A\to\Z$.  Since $A$ is a minimal compatible
set, for any $b\in A$ there is a finite sequence $a_1,\dots,a_r$ of elements of $A$, where
$a=a_1$ and $b=a_r$, such that (for each $j$) $a_j$ and $a_{j+1}$ have the same image under
either $f_1$ or $f_2$.
If $f_i(a_j)=f_i(a_{j+1})$ then $\ell(a_j)/m_i(a_j) = \ell(a_{j+1})/m_i(a_{j+1})$,
so $\ell(a_{j+1}) = \ell(a_j) m_i(a_{j+1})/m_i(a_j)$.  Thus, we can express $\ell(b)$ as
$\ell(a)$ times
a rational number whose numerator and denominator are products of values of $m_1$ and $m_2$.
It follows that any other compatible labeling must be a rational number times $\ell$.  Conversely,
a rational multiple of $\ell$ is a consistent labeling if and only $\ell(b)/m_i(b)\in\Z$ for every
$b\in A$ and $i\in\{1,2\}$.  The result follows.
\end{proof}

\begin{corollary}
\label{zwb}
Suppose $f_1$ and $f_2$ have a common composite, and let $h$ be a common composite of minimal
degree.  Then the minimal compatible sets $A\subset\Kbar$ are precisely the sets
$\{b\in\Kbar:h(b)=h(a)\}$ with $a\in\Kbar$.
Moreover, if $\ell_0$ is the minimal consistent labeling on $A$, then
$\ell_0(a)$ is the multiplicity of $x=a$ as a root of $h(x)-h(a)$, and furthermore
$\sum_{a\in A}\ell_0(a)=\deg(h)$.  Finally, writing
$\hat h:=\prod_{a\in A}(x-a)^{\ell_0(a)}$, we have $\hat h(x)-\hat h(0)\in K[x]$,
and there is a degree-one $\mu\in K[x]$ such that $\hat h(x) = \hat h(0) + \mu(h(x))$.
\end{corollary}
\begin{proof}
For any $a\in\Kbar$, let $\ell(a)$ denote the ramification index of $x=a$ in
$\Kbar(x)/\Kbar(h(x))$.  The fiber $S=\{b\in\Kbar:h(b)=h(a)\}$ is compatible, and
$\ell$ is a consistent labeling on $S$.
Note that $\sum_{b\in S} \ell(b) = \deg(h)$.  Let $A$ be a minimal
compatible set contained in $S$.  Then Theorem~\ref{comp} implies that
$\hat h:=\prod_{b\in A}(x-b)^{\ell(b)}$ is a common composite over $\Kbar$.  By minimality of
$\deg(h)$, we must have $\deg(h)\le\deg(\hat h)$, so $A=S$ and $\deg(h)=\deg(\hat h)$.
Likewise, $\ell$ must be the minimal consistent labeling on $A$, since otherwise using a smaller
labeling in Theorem~\ref{comp} would produce a common composite of degree lower than $\deg(h)$.
Now $\hat h(x) = \hat h(0) + \mu(h(x))$ for some degree-one $\mu\in\Kbar[x]$.  Since
$\hat h$ is monic and $h\in K[x]$, the leading coefficient of $\mu$ must be in $K$.
Since the constant terms of both $h$ and $(\hat h(x)-\hat h(0))$ are in $K$, we have
$\mu(0)\in K$.  This completes the proof.
\end{proof}


\section{Inconsistent Sets}
\label{sec inc}

In this section we give examples of $f_1,f_2\in K[x]$ for which there is an
inconsistent subset of $\Kbar$.  By Proposition~\ref{conscomp}, this implies there is
no common composite.  We begin by reworking Example~\ref{strange}.

\begin{exampleme}
\label{firstexample}
Consider $f_1=x^2-x$ and $f_2=x^3-x^2$ over any field $K$.
We claim that $\{0,1\}$ is inconsistent.  For, suppose there were a function $\ell$ on
$\{0,1\}$ satisfying the properties of Definition \ref{cons}.  Since
$f_1(0)=f_1(1)$ and $f_2(0)=f_2(1)$, we would have
$$
\frac{\ell(0)}{m_1(0)}=\frac{\ell(1)}{m_1(1)}\qquad\text{ and }\qquad
\frac{\ell(0)}{m_2(0)}=\frac{\ell(1)}{m_2(1)},
$$
so
$$
\frac{m_1(0)}{m_1(1)} = \frac{\ell(0)}{\ell(1)} = \frac{m_2(0)}{m_2(1)}.
$$
But $m_1(0)=m_1(1)=m_2(1)=1$ and $m_2(0)=2$, contradiction.
\end{exampleme}

In the above example the set $\{0,1\}$ is compatible, but this property
is not used in proving there is no common composite.  (By contrast, we crucially used
this property when we treated these polynomials in Example~\ref{strange}.)
It is not difficult to construct
similar examples involving noncompatible inconsistent sets---for instance, one could
replace $f_1$ by $(x^2-x)(x^2-x-1)$.

Our next example involves a larger inconsistent set.

\begin{exampleme}
Consider $f_1=x^3+x+1$ and $f_2=x^4+x+1$ in $\F_3[x]$.
We claim that $A:=\{0,-1,i,i-1\}$ is inconsistent.  For, suppose there is a
consistent labeling $\ell$ on $A$.
Since $f_1(i)=f_1(0)=1$ and $m_1(i)=m_1(0)=1$, we have $\ell(i)=\ell(0)$.
Since $f_2(0)=f_2(-1)=1$ and $m_2(0)=1$ and $m_2(-1)=3$, we have
$\ell(-1)=3\ell(0)$.  Since $f_1(-1)=f_1(i-1)=-1$ and $m_1(-1)=m_1(i-1)=1$,
we have $\ell(i-1)=\ell(-1)$.  Since $f_2(i-1)=f_2(i)=i-1$ and $m_2(i-1)=m_2(i)=1$,
we have $\ell(i)=\ell(i-1)$.  Thus
$$
\ell(i)=\ell(i-1)=\ell(-1)=3\ell(0)=3\ell(i),
$$
contradicting the fact that $\ell(i)$ is nonzero.
\end{exampleme}

These two examples generalize as follows:
\begin{theorem}
\label{cycleram}
Suppose $c_1,\dots,c_{2d}\in\Kbar$ satisfy $f_1(c_i)=f_1(c_{i+1})$ for odd $i$ and
$f_2(c_i)=f_2(c_{i+1})$ for even $i$ (where $c_{2d+1}:=c_1$).
If $f_1$ and $f_2$ have a common composite then
\begin{equation}
\label{cyc}
1 = \prod_{i=1}^d  \frac{m_1(c_{2i-1})}{m_2(c_{2i-1})} \frac{m_2(c_{2i})}{m_1(c_{2i})}.
\end{equation}
\end{theorem}

\begin{proof}
Suppose $f_1$ and $f_2$ have a common composite, and let $\ell$ be a consistent labeling
on $\Kbar$.  Then $\ell(c_i)/\ell(c_{i+1})$ equals $m_1(c_i)/m_1(c_{i+1})$ if $i$ odd,
and equals $m_2(c_i)/m_2(c_{i+1})$ otherwise.  The desired formula follows by
computing the product of all $2d$ terms $\ell(c_i)/\ell(c_{i+1})$.
\end{proof}

We do not know how often one can satisfy the criteria of this Proposition.
Namely, if one begins with a value $c_1$ such that $m_1(c_1)>1$ (i.e., $f_1'(c_1)=0$),
then how likely is it that there exist $c_2,\dots,c_{2d}$ such that
$f_1(c_i)=f_1(c_{i+1})$ for odd $i$ and
$f_2(c_i)=f_2(c_{i+1})$ for even $i$?  If such $c_i$ do exist, one would expect
that `usually' Equation~(\ref{cyc}) is not satisfied.  However, we suspect that
it is rare for such $c_i$ to exist.

As an extreme example in this direction, we note that there are polynomials $f_i$
for which $\Kbar$ is consistent, even though the $f_i$ have no common composite:

\begin{exampleme}
Let $f_1=x^2$ and $f_2=(x-1)^2$ be polynomials over\/ $\Q$.  
Then $m_i(\alpha)=1$ for all $\alpha\in\Qbar$ and $i\in\{1,2\}$,
except that $m_1(0)=2$ and $m_2(1)=2$.  Thus, the constant function $\ell=2$ is a
consistent labeling on\/ $\Qbar$.  However, any compatible subset $S$ of $\Qbar$ would
have to be closed under the map $x\mapsto -x$ (since $-x$ and $x$ are in the same
fiber of $f_1$), and likewise $S$ would be closed under $x\mapsto 2-x$.  But then
$S$ would be closed under the composite map $x\mapsto 2+x$, contradicting finiteness
of $S$.  Hence there is no compatible subset of\/ $\Qbar$, so $f_1$ and $f_2$ have
no common composite.
\end{exampleme}


\section{Derivatives}
\label{sec der}

In the previous section we gave a method which, in certain special cases,
enables one to prove that two polynomials $f_1$ and $f_2$ have no common composite.
In this section we give a more robust method for this.

\begin{proposition}
\label{der}
Suppose that $f_1,f_2\in K[x]$ have a common composite $h$, and suppose
$\alpha,\beta\in\Kbar$ satisfy $h'(\alpha)h'(\beta)\ne 0$ and
$f_i(\alpha)=f_i(\beta)$ for both $i=1$ and $i=2$.
Then $f_1'(\alpha)f_2'(\beta) = f_1'(\beta)f_2'(\alpha)$.
\end{proposition}

\begin{proof}
Writing $h=F_i\circ f_i$ with $F_i\in K[x]$, we have
\begin{align*}
h'(\alpha) &= F_i'(f_i(\alpha))\cdot f_i'(\alpha) \\
h'(\beta) &= F_i'(f_i(\beta))\cdot f_i'(\beta) = F_i'(f_i(\alpha))\cdot f_i'(\beta).
\end{align*}
Since $h'(\beta)\ne 0$, this implies
\[
\frac{h'(\alpha)}{h'(\beta)} = \frac{f_i'(\alpha)}{f_i'(\beta)}.
\]
Since the left side of this equation does not depend on $i$, the result follows.
\end{proof}

\begin{exampleme}
Consider $f_1=x^3$ and $f_2=x^2+x$ over $K=\F_2$.  Letting $\omega$ be a primitive
cube root of unity in $\Kbar$, we see that $f_i(\omega^j)=1$ for each $i,j\in\{1,2\}$.
Since $f_1'(\omega)f_2'(\omega^2) \ne f_1'(\omega^2)f_2'(\omega$),
Proposition~\ref{der} implies that every common composite $h$ of $f_1$ and $f_2$ must
satisfy $h'(\omega)h'(\omega^2)=0$.  In this instance, we know by Proposition~\ref{funny}
that $f_1$ and $f_2$ have a common composite, and that a minimal-degree common composite
is $\hat h := (x^4+x)^3$.  And indeed, $\hat h'(\omega)=\hat h'(\omega^2)=0$.
\end{exampleme}

This example illustrates how to use Proposition~\ref{der} to prove a property of
common composites, assuming such composites exist.  We now build this into a
criterion enabling us to prove nonexistence of a common composite in some cases.

\begin{lemma}
Suppose $f_1,f_2\in K[x]\setminus K[x^p]$ have a common composite, and let $h$ be a
minimal-degree common composite.
For any $\alpha\in\Kbar$ such that $[K(\alpha)\col K]$ is divisible by a prime greater than
\,$\max(\deg(f_1),\deg(f_2))$, we have $h'(\alpha)\ne 0$.
\end{lemma}

\begin{proof}
By Proposition~\ref{conscomp}, there is a compatible consistent set $A\subset\Kbar$
containing $\alpha$.  Assume $A$ is the minimal such set; then $A$ consists of all
$\beta\in\Kbar$ for which there is a finite sequence of elements of $\Kbar$, starting
with $\alpha$ and ending with $\beta$, such that consecutive members of the sequence
have the same image under either $f_1$ or $f_2$.  Our condition on the degrees implies
that the large prime dividing $[K(\alpha)\col K]$ also divides $[K(\gamma)\col K]$ for each
$\gamma$ in the sequence, so this prime divides $[K(\beta)\col K]$, whence $f_i'(\beta)\ne 0$.
Thus $\ell=1$ is the minimal consistent labeling on $A$, so Corollary~\ref{zwb} implies
that $\hat h(x):=\prod_{a\in A}(x-a)$ satisfies $\hat h(x)-\hat h(0)=\mu(h(x))$
for some degree-one $\mu\in K[x]$.  In particular, since $\hat h'(\alpha)\ne 0$, we must have
$h'(\alpha)\ne 0$.
\end{proof}

Combining the previous two results gives our desired criterion:

\begin{corollary}
\label{der2}
Suppose $f_1,f_2\in K[x]\setminus K[x^p]$ and $\alpha,\beta\in\Kbar$ satisfy
$f_i(\alpha)=f_i(\beta)$ for both $i=1$ and $i=2$, and also
$[K(\alpha)\col K]$ is divisible by a prime greater than
\,$\max(\deg(f_1),\deg(f_2))$.  If
$f_1'(\alpha)f_2'(\beta)\ne f_1'(\beta)f_2'(\alpha)$
then $f_1$ and $f_2$ have no common composite.
\end{corollary}

\begin{exampleme}
Consider $f_1=x^4+x^3$ and $f_2=x^6+x^2+x$ over $\F_2$.  One can check that
$\psi(x):= x^{14} + x^{10} + x^9 + x^8 + x^7 + x^6 + x^4 + x + 1$ is irreducible
over $\F_2$.  For any root $\alpha$ of $\psi$, let $\beta=\alpha^{128}$.  Then
$f_i(\alpha)=f_i(\beta)$ for each $i$,
but $f_1'(\alpha)f_2'(\beta)\ne f_1'(\beta)f_2'(\alpha)$,
so Corollary~\ref{der2} implies the $f_i$ have no common composite.
\end{exampleme}

Our proof of Corollary~\ref{der2} generalizes at once to prove the following:

\begin{theorem}
\label{der3}
For $f_1,f_2\in K[x]\setminus K[x^p]$, suppose $c_1,\dots,c_{2d}\in\Kbar$ satisfy
$f_1(c_i)=f_1(c_{i+1})$ for odd $i$ and
$f_2(c_i)=f_2(c_{i+1})$ for even $i$ (where we define $c_{2d+1}:=c_1$).
Suppose further that $[K(c_1)\col K]$ is divisible by a prime greater than
\,$\max(\deg(f_1),\deg(f_2))$.  If
\[
\prod_{i=1}^d \left(f_1'(c_{2i-1}) f_2'(c_{2i})\right) \,\ne\,
 \prod_{i=1}^d \left(f_2'(c_{2i-1}) f_1'(c_{2i})\right)
\]
then $f_1$ and $f_2$ have no common composite.
\end{theorem}

One can check that there is no loss in only applying this result when the $c_i$ are distinct.

\begin{exampleme}
Consider $f_1=x^2+x$ and $f_2=x^4+x^3+x$ over $\F_2$.
The two primitive cube roots of unity have the same image as one another under
both $f_1$ and $f_2$, but they have degree $2$ over $\F_2$ so the above result
does not apply.  For $d<5$, this is the only choice of distinct $c_i$'s such that
$f_1(c_i)=f_1(c_{i+1})$ for odd $i$ and $f_2(c_i)=f_2(c_{i+1})$ for even $i$.
But for $d=5$ we can choose
$(c_1,\dots,c_{10}):=( w, w^{268}, w^4, w^{49}, w^{16}, w^{196}, w^{64}, w^{784},
w^{256}, w^{67})$  where  $w^{10} + w^9 + w^4 + w^2 = 1$.
Since $[\F_2(w)\col \F_2]=10$, these $c_i$ satisfy all the hypotheses of Theorem~\ref{der3},
so $f_1$ and $f_2$ have no common composite.
\end{exampleme}

We suspect that Theorem~\ref{der3} applies to `most' pairs of
polynomials over a finite field.  This intuition has been reinforced by various
examples we have computed.  Our intuition is based on the following
reasoning: the $c_i$ are defined by $2d$ equations in $2d$ variables,
so `at random' we expect to find solutions.  Specifically, we can
apply the fiber-finding algorithm to the indeterminate $\alpha=t$ in
$K[t]$.  This gives a polynomial $r_{2d}\in K[t,x]$ such that
$r_{2d}(c_1,c_1)=0$, narrowing the choices for $c_1$ to a finite set.
It may happen that no such choice for $c_1$ leads to a solution for
$c_2,\ldots,c_{2d}$ with the $c_i$ distinct, but this seems unlikely 
to happen except in unusual circumstances.  Finally, as we vary $d$, it seems there
should be some $d$ for which a corresponding $c_1$ is defined over an extension of $K$
of degree divisible by a large prime, and moreover `at random' the products of derivatives
expressed in Theorem~\ref{der3} are almost certainly distinct.

Unfortunately, there are cases where two polynomials have no common composite, but this
nonexistence cannot be proved with Theorem~\ref{der3}.

\begin{exampleme}
Consider $f_1=x^2+x$ and $f_2=x^6+x$ over $\F_2$.  Since $f_i'(x)=1$, there are no $c_j$'s
satisfying the hypotheses of Theorem~\ref{der3}.
\end{exampleme}

In a subsequent paper we will develop further methods for proving nonexistence of a common
composite, and in particular we will show that the polynomials in the above example have
no common composite.

\end{document}